\UseRawInputEncoding
\documentclass[conference]{IEEEtran}
\IEEEoverridecommandlockouts
\usepackage{cite}
\usepackage{amsmath,amssymb,amsfonts}
\usepackage{algorithmic}
\usepackage{graphicx}
\usepackage{textcomp}
\usepackage{xcolor}
\usepackage[ruled,vlined]{algorithm2e}
\def\BibTeX{{\rm B\kern-.05em{\sc i\kern-.025em b}\kern-.08em
    T\kern-.1667em\lower.7ex\hbox{E}\kern-.125emX}}
\begin{document}
\title{Convex Optimization in Legged Robots}
\author{\IEEEauthorblockN{Prathamesh Saraf$^{1}$, Mustafa Shaikh$^{1}$, Myron Phan$^{2}$}
\IEEEauthorblockA{\textit{$^{1}$Department of Electrical and Computer Engineering} \\
\textit{$^{2}$Department of Mechanical and Aerospace Engineering} \\
\textit{University of California - San Diego}\\
mphan@ucsd.edu, mushaikh@ucsd.edu, psaraf@ucsd.edu}
* All authors have equal contribution}
\maketitle
\begin{abstract}
Convex optimization is crucial in controlling legged robots, where stability and optimal control are vital. Many control problems can be formulated as convex optimization problems, with a convex cost function and constraints capturing system dynamics. Our review focuses on active balancing problems and presents a general framework for formulating them as second-order cone programming (SOCP) for robustness and efficiency with existing interior point algorithms. We then discuss some prior work around the Zero Moment Point stability criterion, Linear Quadratic Regulator Control, and then the feedback model predictive control (MPC) approach to improve prediction accuracy and reduce computational costs. Finally, these techniques are applied to stabilize the robot for jumping and landing tasks.\\
\indent Further research in convex optimization of legged robots can have a significant societal impact. It can lead to improved gait planning and active balancing which enhances their ability to navigate complex environments, assist in search and rescue operations and perform tasks in hazardous environments. These advancements have the potential to revolutionize industries and help humans in daily life.
\end{abstract}
\begin{IEEEkeywords}
convex optimization, legged robots, model predictive control, stability
\end{IEEEkeywords}
\section{Introduction}
Control problems can be formulated as optimization problems by defining an objective function that quantifies the desired behavior of the system, and a set of constraints that capture the physical limitations of the system and any other relevant constraints. The objective function can be some sort of performance measure, such as minimizing energy consumption or maintaining the stability of the system. The constraints include the dynamics of the system, input/output constraints, and state constraints. Once the control problem is formulated as an optimization problem, the goal is to find the inputs or control actions that optimize the objective function while satisfying the constraints. Well-known optimization techniques can be used to solve the resulting optimization problem and find the optimal inputs or control actions that achieve the desired behavior of the system. This approach enables the design of controllers that can handle complex and nonlinear systems and can provide improved performance and robustness compared to traditional control design approaches. Attaining efficient and robust solutions for the optimal control sequence is key in legged robots as the computation power available is limited, and trajectory generation is only one of many tasks the on-board computer must complete. Therefore, by formulating the control problems as convex problems, designers can take advantage of highly developed and efficient solvers to obtain the optimal solution.

\section{Convex Optimization Applications}
\label{prob-formulation}
We start with some literature and initial works on convex optimization applications in legged robots, which lay the foundation to the most widely used optimization methods, Model Predictive Control.
\subsection{Zero Moment Point}
\label{zmp}
\begin{figure}
    \centering
    \includegraphics[height=0.3\textwidth,width=0.4\textwidth]{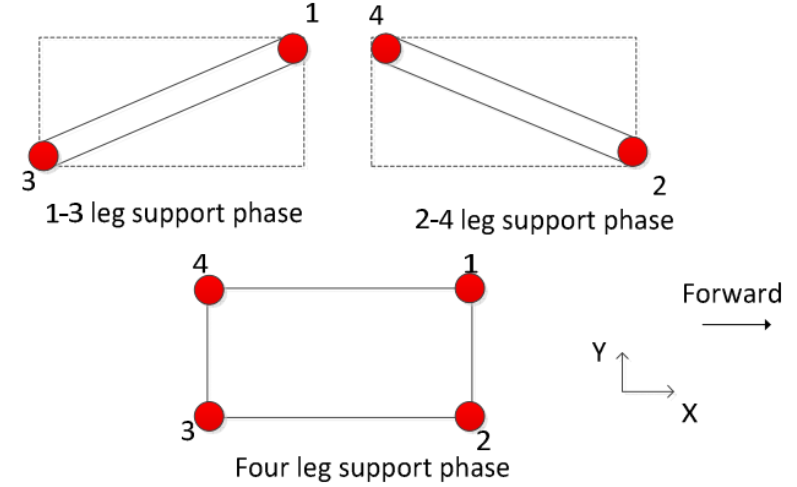}
    \caption{Support Polygon for ZMP}
    \label{fig:my_label}
\end{figure}

We first briefly introduce some background on the necessary dynamics needed to understand legged robot motion. The Zero Moment Point (ZMP) criterion offers several advantages for analyzing the stability of legged robots. It is a point/projection on the ground at which all the moments and the forces acting on the robot sum up to zero. It is a simple and intuitive method that is easy to understand and apply. It enables real-time monitoring, allowing legged robots to adapt to changes in their environment and maintain balance in dynamic scenarios\cite{b13}. Additionally, the ZMP criterion can be combined with whole-body control optimization algorithms to compute optimal joint angles and control inputs, ensuring stability across different types of legged robots \cite{b14}. Figure 1 shows the support polygon cases and Figure 2 illustrates the ZMP criterion and friction cone constraints in simulation with respect to the ANYmal quadruped robot. The ZMP stability optimization equations are given below:
\begin{equation}
\min \| (\boldsymbol{x}_{\mathrm{zmp}} - \boldsymbol{x}_{\mathrm{cp}}),( \boldsymbol{y}_{\mathrm{zmp}} - \boldsymbol{y}_{\mathrm{cp}}) \|
\end{equation}
subject to
\begin{equation}
z \ddot{x}-\left(x-x_{z m p}\right)(\ddot{z}+g) = 0
\end{equation}

where,
\begin{equation}
\begin{aligned}
X_{z m p} & =X \mathrm{com}-\frac{Z \mathrm{com}}{g} \ddot{X c o m} \\
Y_{z m p} & =Y \mathrm{com}-\frac{Z \mathrm{com}}{g} \ddot{Y} \mathrm{com}
\end{aligned}
\end{equation}
However, the ZMP criterion has certain limitations. It assumes a simplified model of a rigid body with a fixed center of mass on a flat and rigid surface, which may not accurately represent real-world conditions. Legged robots often encounter uneven or slippery surfaces and may have flexible components that affect their balance. Moreover, in dynamic environments, such as crowded streets or rough terrain, the ZMP criterion may not account for unexpected obstacles or external forces that can impact stability. Additionally, the ZMP algorithm relies on precise measurements of position, velocity, dynamics, and control inputs, making it sensitive to sensor noise, communication delays, and other sources of uncertainty.
\begin{figure}
    \centering
    \includegraphics[height=0.3\textwidth,width=0.4\textwidth]{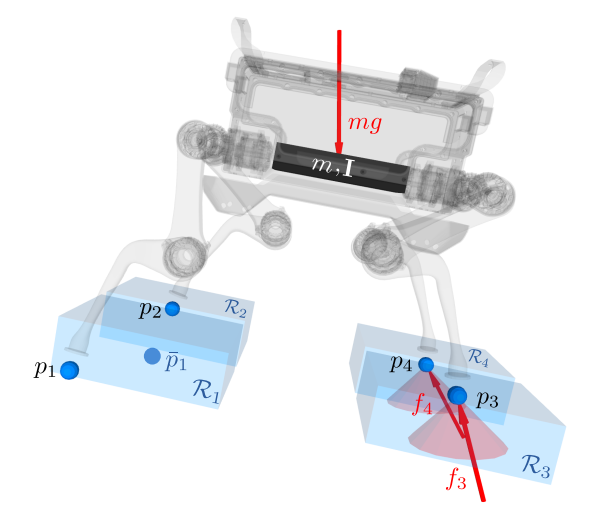}
    \caption{ZMP visualization on a quadruped robot}
    \label{fig:my_label}
\end{figure}
While the ZMP criterion provides valuable insights into the stability of legged robots, its limitations in considering real-world complexities and uncertainties emphasize the need for more advanced algorithms and techniques to ensure robust and reliable stability control in legged robot applications. A more improved version of the ZMP problem gives rise to the Linear Quadratic Regulator control which is described in the next section.

\subsection{Linear Quadratic Regulator}
\label{lqr}
LQR control, or Linear Quadratic Regulator control, is a popular control strategy used in legged robots to achieve stability. By optimizing a quadratic cost function, LQR control computes control inputs that minimize the deviation from desired states, ensuring stability and smooth movements. LQR control is known for its simplicity and ease of implementation, making it a widely adopted approach in legged robot research. The optimization cost function is given below:
\begin{equation}
\begin{gathered}
J=\int_0^{\infty}\left(x^T Q x+u^T R u\right) d t \\
A^T P+P A-P B R^{-1} B^T P+Q=0 \\
K=R^{-1} B^T P \\
\tau=\tau_0-K x .
\end{gathered}
\end{equation}
where J is the cost function to be optimized, Q and R are PSD weight matrices and x is the state of the system with u being the system input. $\tau$ represents the torque commands required for leg joints.
One advantage in terms of legged robots is their robustness, as LQR control can be designed to handle disturbances and uncertainties in the environment, allowing legged robots to maintain stability even in the presence of unexpected changes\cite{b11}. Additionally, LQR control is relatively easy to implement and tune, making it a popular choice for researchers and engineers working on stability control in legged robots\cite{b12}.

However, LQR control is limited to linear systems \cite{b24} \cite{b25}, which can restrict its effectiveness in complex and nonlinear systems commonly encountered in legged robots. Another drawback is the requirement for an accurate model of the legged robot's dynamics for LQR control to be effective, which can be challenging to obtain. Furthermore, LQR control can be computationally intensive for large and complex-legged robots, leading to reduced real-time performance and limited adaptability to sudden changes in the environment.

Though LQR control offers advantages such as robustness and ease of implementation, its limitations in handling nonlinear systems, reliance on accurate models, and computational complexity highlight the need for more advanced control strategies for stability control in legged robots. We thus move on to the Model Predictive control technique which is widely used in current stability algorithms.

\subsection{Model Predictive Control}
\label{prob-formulation}
In the paper \cite{b2}, Jared Di Carlo, et al. explain how the control problem for the MIT Cheetah 3 robot is formulated as an optimization problem using model-predictive control (MPC). The objective of the MPC is to minimize a cost function that captures the desired behavior of the robot, such as maintaining stability and achieving high-speed locomotion. The cost function is subject to a set of constraints that capture the dynamics of the robot, control constraints, and state constraints. The MPC problem is solved using a linear-quadratic program (LQP) formulation, which is a type of convex optimization problem that involves minimizing a quadratic objective function subject to linear constraints. The LQP formulation incorporates a simplified dynamic model of the MIT Cheetah 3, which includes a set of nonlinear constraints that capture the physical limitations of the robot. The authors use a convex approximation of these nonlinear constraints, which allows them to formulate the MPC problem as an LQP that can be efficiently solved using standard solvers. This approach enables the design of controllers that can achieve stable and dynamic locomotion in the MIT Cheetah 3 robot, even under challenging conditions such as uneven terrain and disturbances.

 The objective of the LQP is to minimize a cost function that captures the desired behavior of the robot, such as maintaining stability and achieving high-speed locomotion. The cost function is subject to a set of linear constraints that capture the dynamics of the robot, control constraints, and state constraints. The authors demonstrate the effectiveness of the convex MPC approach by testing it on the MIT Cheetah 3 robot, which is a highly dynamic quadrupedal robot capable of high-speed locomotion and agile maneuvers. The results show that the convex MPC approach is able to achieve stable and dynamic locomotion in the robot, even under challenging conditions such as uneven terrain and disturbances.

In this paper \cite{b3}, the authors formulate the MPC problem as a quadratic program, which is a convex optimization problem that can be efficiently solved using standard solvers. The objective of the quadratic program is to minimize a cost function subject to a set of constraints, which include the dynamics of the system, control constraints, and state constraints. The authors use a simplified linear model of the system dynamics for prediction, which allows them to formulate the MPC problem as a quadratic program that can be efficiently solved. However, the accuracy of the predictions is improved by incorporating feedback from a low-level controller, which corrects the predicted trajectories and reduces the computational complexity of the MPC. The use of convex optimization techniques, such as quadratic programming, enables the authors to solve the MPC problem efficiently and effectively and to achieve better control performance compared to traditional MPC approaches. The use of convex optimization in this paper is an important contribution to the field of legged robot control and demonstrates the power of optimization techniques in addressing challenging control problems in robotics. The MPC formulation goes as below:
\begin{equation}
\begin{array}{rc}
\min _{\mathbf{x}, \mathbf{u}} & \sum_{i=0}^{k-1}\left\|\mathbf{x}_{i+1}-\mathbf{x}_{i+1, \mathrm{ref}}\right\|_{\mathbf{Q}_i}+\left\|\mathbf{u}_i\right\|_{\mathbf{R}_i} \\
\end{array}
\end{equation}
\begin{equation}
\begin{array}{rc}
\text { subject to } \quad \mathbf{x}_{i+1}=\mathbf{A}_i \mathbf{x}_i+\mathbf{B}_i \mathbf{u}_i, i=0 \ldots k-1 \\
\mathbf{c}_i \leq \mathbf{C}_i \mathbf{u}_i \leq \overline{\mathbf{c}}_i, i=0 \ldots k-1 \\
\mathbf{D}_i \mathbf{u}_i=0, i=0 \ldots k-1 \\
f_{\min } \leq f_z \leq f_{\max } \\
-\mu f_z \leq f_x \leq \mu f_z \\
-\mu f_z \leq f_y \leq \mu f_z
\end{array}
\end{equation}
The constraints describe the system dynamics and the friction cone constraints for each foot in contact with the ground. The optimization equation minimizes the error in the current state and the reference trajectory which is computed offline initially. Q and R are positive semi-definite weight matrices.

Thus we see that Model Predictive Control (MPC) offers several advantages for stability control in legged robots \cite{b23}. First, it enables optimized control by minimizing a cost function that represents stability. This allows for the generation of control inputs that ensure the robot maintains stability and achieves desired performance. Additionally, MPC is more robust to disturbances and uncertainties in the environment compared to Linear Quadratic Regulator (LQR) control. It can handle nonlinear systems, providing more flexibility in modeling the complex dynamics of legged robots.

Though, one of the main challenges is its computational complexity, especially for large and complex-legged robots. The computational demands of MPC can affect real-time performance and the robot's ability to quickly respond and adapt to changes in the environment. The look-ahead horizon control techniques can mitigate this challenge by allowing the robot to pre-plan its behavior and optimize control inputs in advance, thus improving the robot's ability to handle real-time control tasks efficiently.

\subsection{Sequential Linear Quadratic - Model Predictive Control}
Sequential Linear Quadratic (SLQ) Model Predictive Control (MPC) \cite{b19} is an advanced control strategy utilized in legged robots for stability control. SLQ-MPC approximates the nonlinear dynamics of the legged robot with a linear model and computes a linear feedback control law at each time step. By optimizing a quadratic cost function over a finite time horizon, SLQ-MPC enables the computation of optimal control inputs to ensure stability. This approach combines the advantages of MPC, such as optimized control and adaptability to changes in the environment, with the computational efficiency of linear approximations. SLQ-MPC is particularly well-suited for legged robots as it addresses the challenges of modeling and controlling the complex dynamics involved in legged locomotion. By leveraging linear approximations, SLQ-MPC provides a practical and efficient solution for stability control in legged robots while still achieving high-performance results. This also laid down the foundation for the use of second-order convex problems which have much faster computation and robustness compared to the classical optimization algorithms.

\subsection{Non-Linear Model Predictive Control}
Non-Linear Model Predictive Control (NMPC) is a sophisticated control strategy employed for the motion control of quadruped robots. NMPC is an optimization-based approach that deals with non-linearities and constraints more effectively than traditional control strategies. In the case of quadruped robots, the system dynamics are non-linear due to the multi-body mechanical structure and its interaction with the environment. NMPC applies an internal model to predict the system's future behavior over a finite prediction horizon, then calculates the control inputs that minimize a specified cost function.

The formulation of NMPC involves determining a cost function that is minimized over the prediction horizon. This cost function represents the discrepancy between the predicted output and the desired output. The non-linear system's dynamics are represented by a set of non-linear differential equations. The control inputs are calculated by solving a non-linear optimization problem at each time step. An important feature of this approach is that it can handle constraints, for example, on the control inputs and states, making it a QCQP problem \cite{b20} \cite{b21}.

The general NMPC formulation can be expressed as follows:

\begin{equation}
\begin{aligned}
& \underset{u(.),x(.)}{\text{min}}
& & \int_{t}^{t+T} L(x(t),u(t))dt + V(x(t+T)) \
\end{aligned}
\end{equation}

\begin{equation}
\begin{aligned}
& \text{subject to}
& & \dot{x}(t) = f(x(t),u(t)), \\
&&& x(t) \in X, ; u(t) \in U, \\
&&& x(t+T) \in X_f, \\
&&& x(t) = x_{0}.
\end{aligned}
\end{equation}

Where $L$ is the Lagrangian (running cost), $V$ is the terminal cost, $f$ is the system dynamics, $x$ are the states, $u$ are the controls, $X$ and $U$ are the state and control constraints respectively, $X_f$ is the set of terminal states, and $x_{0}$ is the initial state.

One significant advantage of NMPC for quadruped robots is its ability to account for the system's non-linear dynamics and constraints in a principled way. This can result in more precise, efficient, and stable robot movement, especially in complex or uncertain environments. It can handle the dynamics and physical constraints of a quadruped robot more effectively than traditional linear control methods, such as PID, LQR, and Linear MPC as it constantly re-evaluates and adapts to new inputs and changes in the environment.

However, the major disadvantage of NMPC is its computational complexity. The need to solve a non-linear optimization problem at each time step can be computationally intensive, particularly for systems with high dimensionality like quadruped robots. This can result in high latency in control commands, making NMPC less suitable for real-time control applications without significant computational resources or simplified approximations. Moreover, the performance of NMPC highly depends on the accurate modeling of the system dynamics, which can be challenging in real-world applications due to uncertainties and disturbances.

\begin{figure}
    \centering
    \includegraphics[height=0.2\textwidth,width=0.4\textwidth]{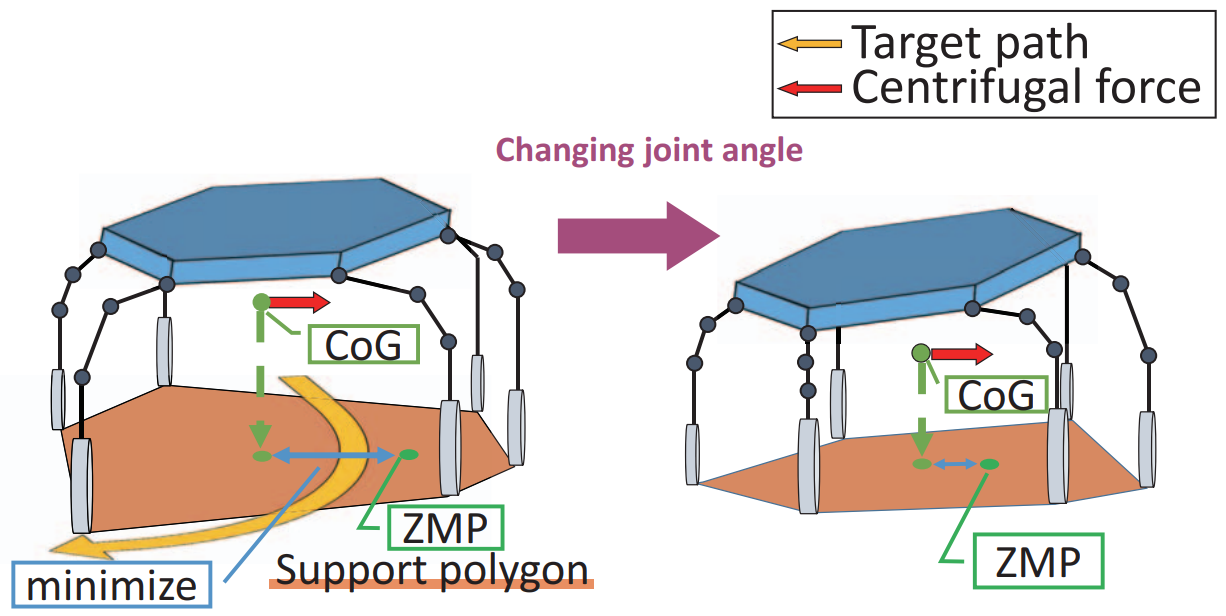}
    \caption{ZMP+MPC approach for posture control - QCQP}
    \label{fig:my_label}
\end{figure}

In \cite{b22}, the authors have formulated the non-linear MPC problem using the ZMP criterion. The Zero Moment Point (ZMP) combined with Model Predictive Control (MPC) provides a powerful formulation for posture correction and trajectory optimization in legged robots \cite{b22}. ZMP is a dynamic stability criterion used extensively in bipedal and quadrupedal robot locomotion, representing the point on the ground where the total moment of the inertial forces and the gravity forces is zero. By maintaining the ZMP within the support polygon of the robot (the area enclosed by the feet in contact with the ground), one can ensure dynamic balance. With MPC, the robot's future behavior is predicted over a finite horizon, and the controls are iteratively updated by solving an optimization problem, ensuring the ZMP stays within the support polygon.

In the ZMP+MPC formulation, the primary objective is to minimize the deviation of the ZMP from a reference trajectory while also considering factors like energy consumption, joint torques, and smoothness of motion. This problem can be posed as a Quadratically Constrained Quadratic Programming (QCQP) problem, where the objective function and constraints are all quadratic. The ZMP constraints ensure the ZMP stays within the support polygon, and the MPC optimization will continually adjust the robot's posture and trajectory to maintain stability, even when the robot is executing dynamic maneuvers or responding to external disturbances. The approach provides a framework for planning dynamically stable and feasible trajectories, which are crucial for successful navigation in complex environments. The equations below depict the non-linear constraints and frame the QCQP optimization problem for MPC stabilization:

\begin{equation}
\begin{aligned}
& J=\phi(\boldsymbol{x}(t+T))+\int_t^{t+T} L(\boldsymbol{x}(\tau), \boldsymbol{u}(\tau)) \mathrm{d} \tau, \\
& \phi(\boldsymbol{x}(t))=Q_{1 \mathrm{f}} \frac{1}{S(\boldsymbol{x}(t))}+Q_{2 \mathrm{f}} e^2(\boldsymbol{x}(t)) \\
& +Q_{3 \mathrm{f}}\left(h_{\mathrm{ref}}-h_{\mathrm{b}}(\boldsymbol{x}(t))\right)^2, \\
& L(\boldsymbol{x}(t), \boldsymbol{u}(t))=Q_1 \frac{1}{S(\boldsymbol{x}(t))}+Q_2 e^2(\boldsymbol{x}(t)) \\
& +Q_3\left(h_{\mathrm{ref}}-h_{\mathrm{b}}(\boldsymbol{x}(t))\right)^2 \\
& +\boldsymbol{u}(t)^T \boldsymbol{R} \boldsymbol{u}(t) \text {. } \\
& e=\sqrt{\left(x_{\mathrm{ZMP}}-x_{\mathrm{CoG}}\right)^2+\left(y_{\mathrm{ZMP}}-y_{\mathrm{CoG}}\right)^2} . \\
&
\end{aligned}
\end{equation}

\subsection{Active Balancing}

We have seen that convex optimization plays an important role in several types of robotics control problems. A particular, general example of this is active balancing problems. In the field of legged robots, it is crucial to ensure the robot remains balanced throughout the course of its trajectory. Active balancing refers to a class of problems that corrects and ensures the stability of the actual robot motion given a trajectory that may not be dynamically feasible or stable. For example, a control input that would cause the robot to fall over should be modified so that the robot remains balanced. This is done by minimizing an objective function that aims to keep the actual trajectory as close as possible to a pre-planned input, while adding constraints that model the dynamics of the problems and keep the robot balanced. In the main paper we will explore in this section,\cite{b1} the authors present a general framework into which many active balancing problems can be framed as second-order cone programming (SOCP) problems.
\begin{figure}
    \centering
    \includegraphics[height=0.1\textwidth,width=0.4\textwidth]{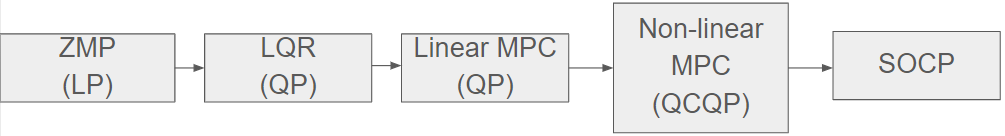}
    \caption{Evolution of optimization algorithms for the stability of legged robots}
    \label{fig:my_label}
\end{figure}
\subsubsection{Previous work}
Previous optimization-based approaches generally use quadratic programs and involve methods such as minimizing least squares tracking errors while maintaining the center of mass at a predefined point. For example, in \cite{b9}, the authors use a quadratic program to optimize the actions a robot takes to maintain stability in the presence of an external impact such as a push or a shove. They optimize an objective function that minimizes the square sum of acceleration of joints, while also ensuring the robot remains standing and does not sit down to maintain stability in the presence of the impact. The objective is given as follows:
\begin{equation}
\text{minimize } \ddot{q}C_q \ddot{q} - s_{\ddot{y}}
\end{equation}
subject to various (generally non-convex) dynamics constraints, and where $\ddot{q}$ represents joint acceleration and $s_{\ddot{y}}$ is the acceleration in the y-direction to ensure the robot remains standing during the optimal solution. Note that in this approach, a limitation is that the authors did not formulate the problem as a convex problem and acknowledge that the optimization has to be done locally, and therefore a global solution is not guaranteed. To build on this work, and to obtain guarantees about global optimality, the aim is to formulate such balancing problems as convex problems and to solve them using known algorithms. The main paper we consider formulates this problem and many others as a convex problem to allow for globally optimal and efficient solutions.
\subsubsection{Current work}
In \cite{b1}, the authors show that many such problems can be formulated as SOCP problems, often by introducing a new variable that helps cast the problem as a SOCP program.

First, we introduce the underlying concepts. A second order cone in $\mathbb{R}^{p+1}$ is defined as $K_p = \{(x,y) \in \mathbb{R}^{p+1} : ||x||_2 \leq y\}$. This set is convex since it is the intersection of an (infinite) number of halfspaces: $K_p = \cap_{u:||u||_2 \leq 1} \{(x,y) \in \mathbb{R}^{p+1}: x^T u \leq y\}$.

\begin{figure}
    \centering
    \includegraphics[height=0.3\textwidth,width=0.2\textwidth]{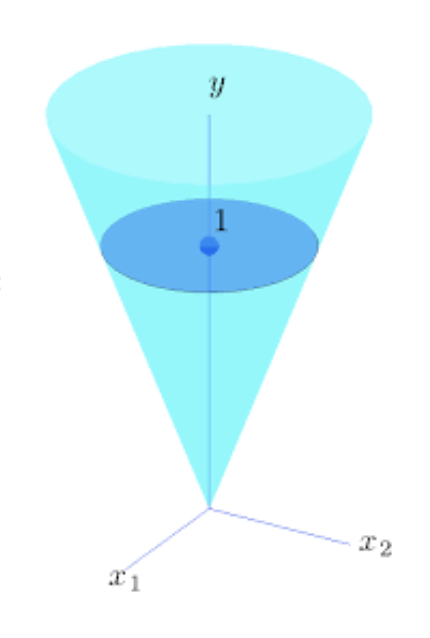}
    \caption{Second Order Cone}
    \label{fig:socp}
\end{figure}

SOCP problems are a class of problems formulated as follows:
\begin{equation}
\begin{aligned}
&\min_{x} f^T x \
\text{subject to} \quad &||A_i x + r_i|| \leq l_i^{T} x + m_i, \quad i = 1,\ldots,N \
&Dx = g
\end{aligned}
\end{equation}

The objective function is convex (linear) and the constraints define a convex set in that the constraints are of the exact form of the second-order cone definition so it is equivalent to requiring the solution to remain within the second-order cone \cite{b1}. Since the second-order cone is convex, so is the constraint set. Due to the generality of SOCPs, many types of convex problems can be reformulated as SOCPs, for example, QCQPs using rotated cones, problems with sums of norms, problems with hyperbolic constraints, and robust least squares (i.e. when there is uncertainty in the data), among others. 

As discussed earlier, the main objective in many active balancing is to minimize the tracking error of the robot's joint accelerations. A reference joint acceleration trajectory $\ddot{q}_{ref}$ is provided, and the actual trajectory is desired to be as close as possible to the reference under the relevant constraints. This leads to an objective of the form:
\begin{equation}
\min_q ||\ddot{q} - \ddot{q}_{ref}||
\end{equation}
To formulate this as a SOCP, we must take some further steps. First, we introduce a dummy variable to minimize and move the objective to a constraint:
\begin{equation}
\begin{aligned}
&\min_{t} \
\text{subject to} \quad &\min_{q} ||\ddot{q} - \ddot{q}_{\text{ref}}|| \leq t
\end{aligned}
\end{equation}
and several other dynamics constraints that will be presented next.
The final step is to convert the objective into the form $f^T x$, and so the authors introduce a variable $\bold{x}$ and a selection vector $f = \begin{bmatrix} 1 & 0 & 0 & .... & 0 \end{bmatrix}$, where $\bold{x} = (t,\ddot{q}, F_1,..., F_m,\lambda_1,...,\lambda_m)$. Now $\min f^T x$ gives us our original problem of $\min t$. The final remaining piece is to ensure the constraints are either conic or affine. While a detailed discussion of dynamics is out of the scope of this paper, we demonstrate a representative example, and then state the full constraint set with the final SOCP formulation. The constraint that represents the requirement that the ZMP remain within the support polygon can be expressed as a linear inequality of the form: 
\begin{equation}
    (A^S A_{ZMP} + b^S b_{ZMP})F_0 \leq 0
\end{equation}
However, this is generally non-conic. To remedy this, the authors use a change of variable and provide a desired ZMP trajectory and instead keep the magnitude of the difference to within an epsilon provided by the user: $||p_{ZMP} - p^d_{ZMP}|| \leq \epsilon_{ZMP}$. This constraint is now a conic constraint of the form shown above in the formulation of the SOCP problem. The final form of the SOCP formulation is as follows:
\begin{equation}
\begin{aligned}
&\min f^T x \
\text{subject to} \
&||S_{\ddot{q}} x - \ddot{q}{\text{ref}}|| \leq S_t x \
&A_s S{\ddot{q}} x \leq b_s \
&(A^s A_{\text{zmp}} + b^s b_{\text{zmp}})F_0 \leq \epsilon_{\text{zmp}} b_{\text{zmp}} F_0 \
&||P_mz - [p_{\text{zmp}}] P_{xy}) F_0|| \leq S_v x - \epsilon_{\text{zmp}} \mathbf{1}{1\times m} S\lambda x
\end{aligned}
\end{equation}
The constraint set above has been reformulated as conic constraints of the form 
$||A_i x + r_i|| \leq l^{T}_i x + m_i$,
$Dx = g$
\subsubsection{Benefits of formulating as SOCPs}
The benefit of formulating balancing problems as SOCP programs is that several robust and efficient interior point algorithms exist to solve such problems. These programs have been shown to converge in 5-50 iterations regardless of the problem dimension \cite{b6}. Though there is an even more general type of program than SOCP, called a Semi-Definite Program (SDP), solvers for these are not as efficient as those for SOCPs, and so it is desirable to formulate a problem as a SOCP to obtain good generality and known efficient solutions. The authors explore the computational performance of their algorithms and show that they can achieve good performance in various classes of active balancing problems. 

\subsubsection{Limitations and Future Work}
Most work currently tracks and minimizes error between input reference trajectories and actual trajectories of joint angles. However, as computational power increases, these optimization problems can incorporate more variables that can provide more appropriate solutions. A limitation of this paper was that the authors only considered a front kick and side kick motion to assess the performance of their SOCP formulation for active balancing. In the future, more complex trajectories, such as running or jumping, should be considered in order to understand the effectiveness of this formulation under general robot motion. Furthermore, the authors acknowledge the need to improve the computational efficiency of the optimizations. We further discuss these ideas in the next sections.

One other important topic of future research is the handling of non-convex constraints with guarantees of global optimality. For example, convexifying an objective or constraint can often involve relaxing an equality constraint to an inequality, or a non-convex inequality constraint to a convex inequality constraint by introducing a slack variable \cite{b18}. In these situations, the global optimality of the resulting problem is not guaranteed except under certain specific conditions. These conditions have been formulated for certain domains such as aerospace, in which common lower bound constraints on rocket thrust can be relaxed and global optimal solutions still guaranteed under linear system dynamics. In fact, even under non-linear system dynamics, which would lead to a non-linear program, global optimal solutions can be guaranteed if the system dynamics can be approximated by a piecewise affine function. Similar research must be carried out for legged robots - this will enable optimal solutions to be quickly found even under relaxed constraints in non-convex problems.

\subsection{Gait Planning}
\indent Gait planning is a vital part of mobile robots as gaits are the motion patterns by which these robots move. Gaits for robots are similar to how animals or how humans move, with running, jumping, and landing. The papers to be reviewed are focused on optimized jumping approaches on quadrupeds.\\
\begin{figure}
    \centering
    \includegraphics[height=0.3\textwidth,width=0.4\textwidth]{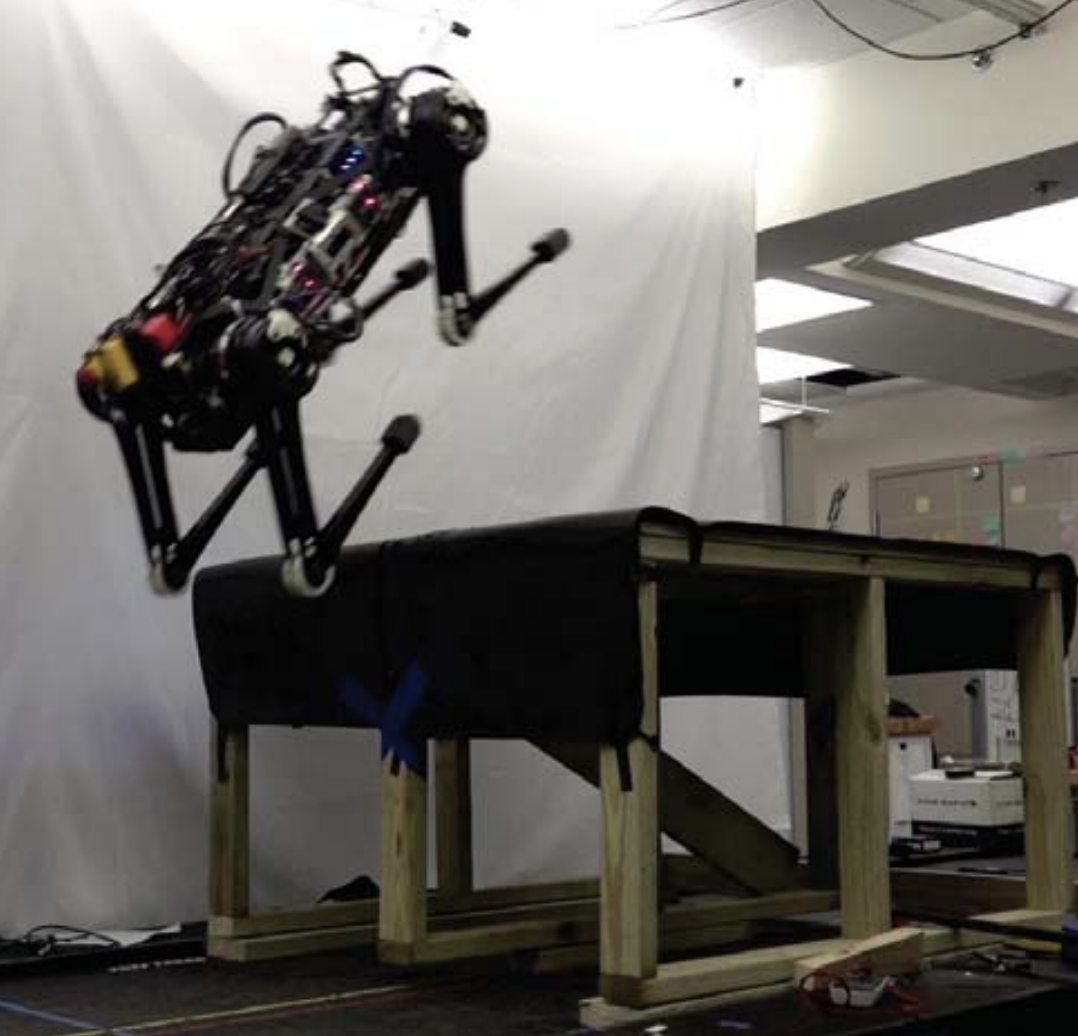}
    \caption{MIT Cheetah 3 Jumping}
    \label{fig:my_label}
\end{figure}
\indent Finding an optimal, real-time approach to gait planning will serve as a step forward toward robust robots in general. Robots that are robust to unknown environments are a crucial functionality for robots to move out of the world of academia towards more real-life applications. There is a particular motivation for gait planning on legged robots because of their inherent ability to navigate difficult terrain and obstacles as opposed to wheeled robots.\\
\indent In general, I think it is hard to understand the importance of this particular subfield of gait planning. Without optimized approaches to planning, we would have to constantly manually tune hyperparameters which will most likely lead to suboptimal results or even failures.

\subsubsection{Online Planning for Autonomous Running Jumps Over Obstacles in High-Speed Quadrupeds}

This paper [15] focused on the real-time optimization problem of jumping. In general, the jumping problem for quadrupeds is a very complicated problem so much of the current research performs offline optimization and stores the results for specific inputs in a controlled test environment. In the context of quadrupeds, this contribution is significant as it was the first to experimentally validate a framework for online planning of running jumps.\\
\indent The dynamics of this problem are the first subject to address. This paper decides to only tackle the problem in 2D as the search space for a 3D problem is often too large for real-time optimization. This means that, theoretically, the quadruped can really only perform well on a flat surface. Most of the testing was done on a flat surface because of this.\\
The equations of motion are:
\begin{align}
    m\Ddot{x}=F_x\\
    m\Ddot{z}=F_z-mg\\
    I\Ddot{\theta}=-xF_x+zF_z
\end{align}
which can be generalized into:
\begin{align}
    \Ddot{x}=u_x\\
    \Ddot{z}=u_z-mg\\
    \Ddot{\theta}=-\alpha xu_x+\alpha zu_z
\end{align}

To make these equations have an analytical solution, we can express these equations in terms of Bezier curves. With these force profiles in terms of Bezier curves, it becomes easy to integrate them to find an analytical trajectory for these equations of motion. With the analytical trajectory, the problem turns into optimizing for the Bezier coefficients to scale the curves.
\begin{align}
    u_x(s)=\sum_{i=0}^n \beta_{i,x}b_{i,n}(s)\\
    u_z(s)=\sum_{i=0}^n \beta_{i,z}b_{i,n}(s)
\end{align}
Once these Bezier coefficients have been optimized, then the analytical solution uses these coefficients (which represent the force profile itself) to generate the actual trajectory.\\
To make the algorithm run in real-time, the optimization foregoes a cost function. Essentially, this becomes a feasibility problem.
\begin{equation}
    \min 0
\end{equation}
with constraints:
\begin{align}
|\theta(t)| &\leq \Theta \
z_{\text{foot}}^f(t) &\leq h_{\text{obs}}(x_{\text{foot}}^f(t)) \
z_{\text{foot}}^h(t) &\leq h_{\text{obs}}(x_{\text{foot}}^h(t)) \
x_{\text{foot}}^f(t_4) &\leq d_0 + 0.5w_0 \
x_{\text{foot}}^h(t_4) &\leq d_0 + 0.5w_0 \
\underline{z} &< z^f(t_1) < \overline{z} \
\underline{z}^h &< z^h(t_2) < \overline{z}^h \
\underline{z} &< z^h(t_3) < \overline{z} \
z^f(t_4) &= \tilde{z}^f \
\underline{\theta} &< \theta(t_4) < \overline{\theta} \
|\alpha_x^j| &< \overline{F}_x \
0 &< \alpha_z^j < \overline{F}_z
\end{align}

\indent The constraints make intuitive sense. Since our state consists of just tracking the center of mass (COM), the leg information is not captured there. To constrain the leg from hitting the obstacle, the first two constraints maintain that the z position of each foot must be above the obstacle. This, of course, is one particular limitation of this algorithm as it relies on knowing prior information about the obstacles themselves.\\
\indent We also want to maintain a certain distance from the obstacle after we land. This is captured by the next two constraints where the x position of each foot must be greater than the distance to the obstacle plus some scaling factor.\\
\indent The next three constraints are related to the physical constraints of the quadruped. The constraints, in words, mean that the position of the legs must be in the workspace of the robot. As in, the position of the legs has to be physically possible for the robot to reach.\\
\indent The next two relate to providing a safe landing configuration for the robot. The final two relate to constraining the ground reaction force that the robot feels on liftoff and landing.\\ 
\indent Overall, this paper provides a good framework for generating real-time trajectories but there are certainly limitations such as the simplified dynamics and fixed force curve profile. While the fixed force profile speeds up convergence, there could potentially be more energy-efficient jump profiles that the robot can execute. A final disadvantage is that the jumps aren't necessarily optimal, in a sense other than the fixed curve profile. Since this is posed as a feasibility problem, the quality of the solution won't be high compared to more complex objective functions. \\
\indent In terms of advantages, the obvious one is that this was the only paper, out of the papers we reviewed, that had an online optimization scheme. There is also another advantage in that this framework allowed the quadruped to make running jumps instead of the more common standing jump. This advantage is also unique and not found in the other papers.\\
\indent In terms of unanswered questions and future works, the authors are unsure of the performance in outdoor terrain obstacles. There are also questions on whether earlier detection of obstacles could lead to better jump results. The authors had implemented a model predictive controller that attempts to place the quadruped in an optimal takeoff position so earlier detection could potentially help.
\begin{figure}
    \centering
    \includegraphics[height=0.3\textwidth,width=0.4\textwidth]{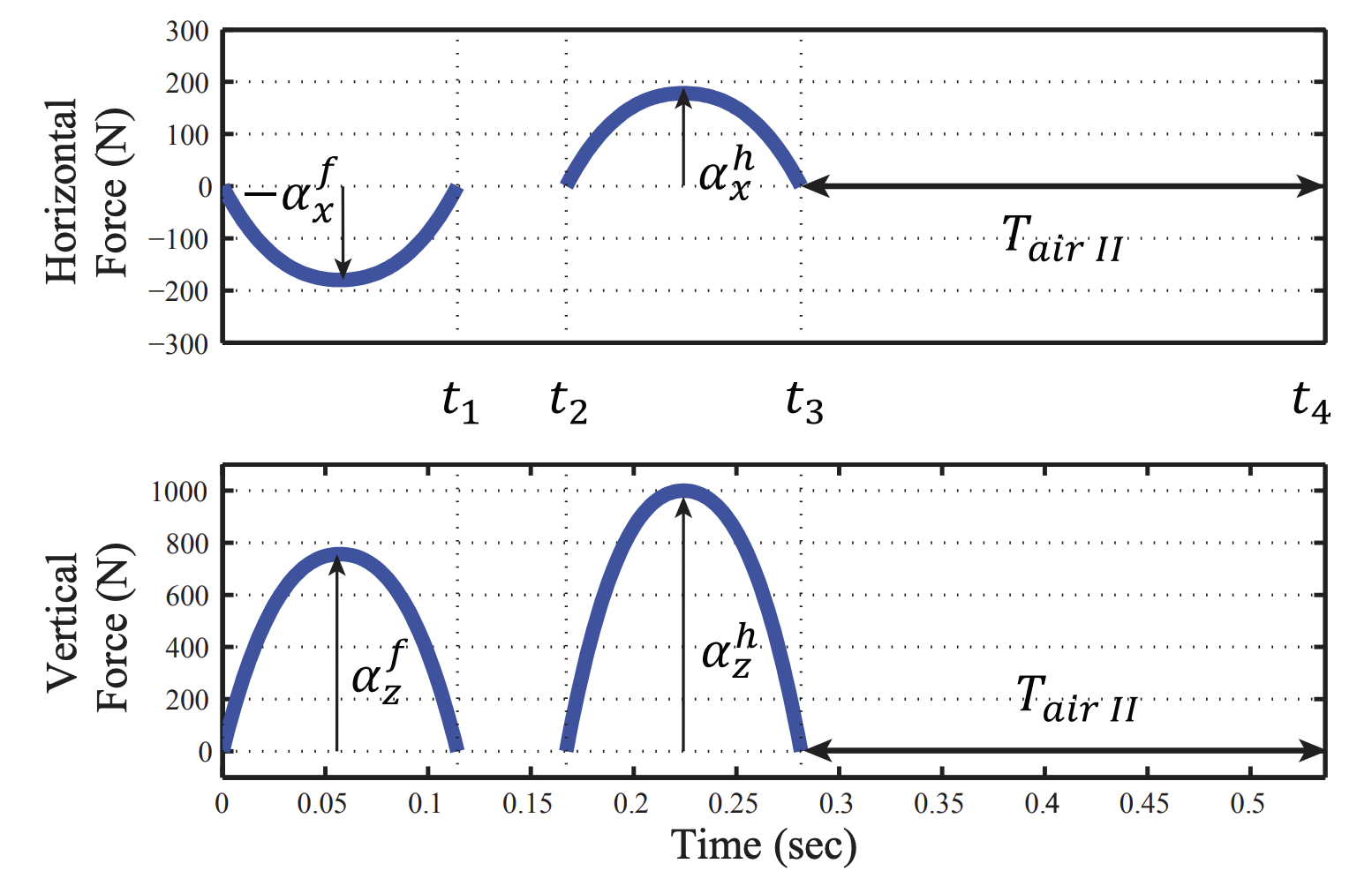}
    \caption{Example Force Profile}
    \label{fig:my_label}
\end{figure}

\subsubsection{Optimized Jumping on the MIT Cheetah 3 Robot}
This paper [16] focuses not just on the feasibility of jumping over an obstacle, it focuses on optimizing the jump itself. This problem is relevant because the problem of there will be many obstacles in outdoor terrain that are much taller than the quadruped itself so it is important to generate an efficient trajectory for high jumping. Out of all the optimized jumping papers that we reviewed, this paper's results had jumped more than 2 times higher than the other results.

Similar to the last paper, the authors chose to constrain this problem in 2D to limit the scope of the problem. The state variables are the 2D position, orientation, and joint angles. The equations of motion are:\\
\begin{align}
    H(Q)\Ddot{Q}+C(Q,\dot{Q})\dot{Q}+g(Q)=\\
    B\tau + B_{fric}\tau_{fric}(\dot{Q})+\sum_i J_i^T(Q)F_i
\end{align}
where H is the mass matrix, C contains the Coriolis and centrifugal terms, B is how torques are entered into the equation, $J_i$ is the spatial Jacobian and $F_i$ is the spatial forces at each foot. This is essentially a generalized force balance where the forces at each foot are mapped to torques.\\
\indent These equations of motion will generate reference trajectories for the optimization problem. Although, the dynamics of the model will change depending on what phase of the jump it is in (all legs on the ground, two legs on the ground, etc.). This, in the previous paper, was not addressed as heavily because performance was not the main concern.\\
\indent Thus, the optimization is split up into three different phases, depending on what phase of the jump it is in. This is because the constraints change. Particularly:
\begin{equation}
    J_{i,stance}(Q)\ddot{Q}+\dot{J}_{i,stance}\dot{Q}=0
\end{equation}
This is the constraint that enforces that a foot is on the ground. This, in words, means that the spatial acceleration of the foot is zero. The actual cost function is a least squares problem:
\begin{multline}
    J=\sum_{k=1}^{N-1}w_q(q_k-q_{ref})^T(q_k-q_{ref}) + w_\tau \tau_k^T\tau_k\\
    +w_N(q_N-q_N^d)^T(q_N-q_N^d)
\end{multline}
With these constraints:
\begin{align}
    q_{min}\le q_k \le q_{max}\\
    |\dot{q}_k|\le\dot{q}_{max}\\
    |\tau_k| \le \tau_{max}\\
    |T_k^x/T^z_k| \le \mu\\
    T^z_k \ge T^z_{min}\\
    q_0=q_{0,d},\dot{q}_0=0\\
    x_0=0, 
    z_0=0,q_{pitch,0}=0,\dot{x}_0=0,\dot{z}_0=0,\dot{q}_{pitch,0}=0\\
    q_k=q_{N,d},\dot{q}_k=0\\
    x_N\ge d_{jumping},\ z_N=h_{platform},q_{pitch,N}=0
\end{align}
\indent These constraints are physical constraints that limit the torque, joint angles, joint velocities, and ground reaction forces.

\indent The first term of the cost function is minimizing the tracking error - how close is the current trajectory from the predicted. The second term seeks to constrain actuator effort. The third term constraints that our final position is where we want it to be.\\
\indent This cost function has some advantages and disadvantages. The advantages are that it allows for the possibility of achieving a higher jump that is optimized for tracking error and actuator effort. Compared to the feasibility problem, this will yield a higher quality solution.\\
\indent A disadvantage is that this algorithm has to be run offline because the optimization is too slow to be run in real-time. There are ways to circumvent this, which the other papers will address, but this paper does not attempt to address them. Another disadvantage is that this paper assumed a planar case, like the last paper. This of course is a disadvantage as the algorithm really only works well on flat surfaces.\\
\indent In terms of unanswered questions / future works, the authors state that they plan to implement vision for autonomous high performance jumping in the future. To do this, a possible unanswered question would be to see if it is possible to modify/implement this algorithm into a real time scenario.\\

\subsubsection{Autonomous Navigation for Quadrupedal Robots with Optimized Jumping through Constrained Obstacles}
This paper [17] aims to perform offline optimization for jumping quadrupeds so that they can jump through a window-shaped obstacle. This contribution is significant/relevant since they seemed to be the first to consider jumping through a constrained obstacle like a window.\\
\indent The dynamics in this paper are largely similar to the previous two papers so I will not review it for this paper. There is some notation to take into account though:
\begin{align}
    q := [q_x,q_z,q_\theta , q_{F1},q_{F2},q_{B1},q_{B2}]^T \in \mathbf{R}^7\\
    x := [q;\dot{q}]\in \mathbf{R}^{14}\\
    u := [\tau_{F1},\tau_{F2},\tau_{B1},\tau_{B2}]^T \in \mathbf{R}^4
\end{align}
\indent Where $q_x,q_z,q_\theta$ represent the planar state and the rest of $\bold{q}$ are the joint angles of the knee and hip. $\bold{x}$ is the state vector that concatenates $\bold{q}$ and its 
derivative. $\bold{u}$ are the torques for the joints represented in $\bold{q}$.\\
The optimization formulation is a little bit different than both of the previous papers. It is still a constrained least squares problem but with additional safety constraints and prediction constraints.
\begin{align}
    \min_{x,\dot{x},u,T} J(x,\dot{x},u,T)\\
    s.t \ x(t_{k+1})=x(t_k)+\frac{\Delta t^{(i}}{2}(\dot{x}(t_{k+1})+\dot{x}(t_k))\\
    \dot{x}(t_{k+1}= f^{(i)}(x(t_k),u(t_k),T(t_k))\\
    x(t_0)=x_0
\end{align}
Where the cost function is defined as:
\begin{multline}
    J=(q(t_{N+1})-q_0)^TP_f(q(t_{N+1})-q_0)+\\ \sum_{k=0}^{N+1}(\dot{q}^T(t_k)Q_{\dot{q}}\dot{q}(t_k)+\ddot{q}^T(t_k)Q_{\ddot{q}}\ddot{q}(t_k))\\+T(t_k)Q_TT(t_k)+u^T(t_k)R_uu(t_k)+\sum_{i=1}^2P_iT_i+\sum_{i=1}^{N_2}P_\delta \delta
\end{multline}
The equality constraints maintain that our solution satisfies our predicted trajectory at the next time step. The first term of the cost function tries to minimize the joint angle error between the takeoff and landing position. This is because the authors defined the starting and end of the jump as a standing position.

\begin{figure}
    \centering
    \includegraphics[height=0.15\textwidth,width=0.4\textwidth]{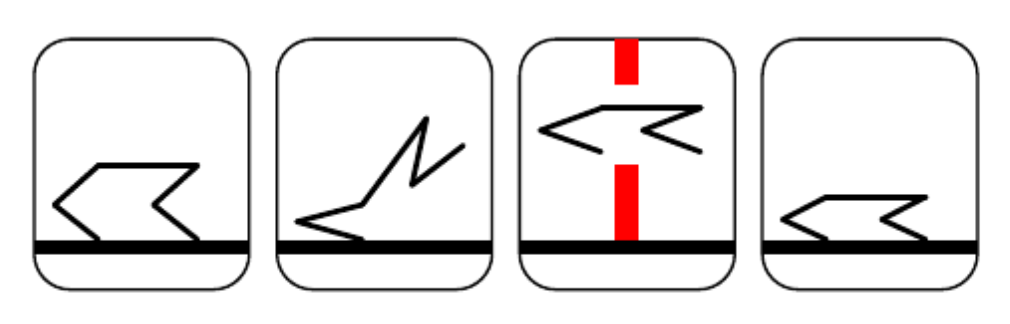}
    \caption{Illustration of jumping through obstacle}
    \label{fig:my_label}
\end{figure}

\indent The second term seeks to minimize joint velocity and acceleration for all time steps. The third and fourth term seeks to minimize ground reaction forces and actuator effort, respectively. The fifth and final terms seek to minimize air time and leg contact, respectively. In short, these terms make sense for the jumping task since minimal joint velocities and accelerations generate smoother trajectories and we always want to constrain actuator effort to some capacity.\\
\indent The authors chose to physically constrain the problem in the usual way of joint angle limits, torque limits, force limits, etc. The constraints are largely the same as the previous papers but more interesting constraints are the ones they chose to prevent the quadruped from jumping into the obstacle ceiling.\\
\indent In essence, the robot itself and obstacles are framed as a convex object if a bounding box is wrapped around it. The robot bounding box must maintain a certain signed distance between itself and any of the obstacles. To achieve this, a new cost function and constraint are added.
\begin{equation}
    J_{w,h}(t_k)=\sum_k w(t_k)^2+h(t_k)^2
\end{equation}
\begin{equation}
    -l^T\mu(t_k)+(A_OP(t_k)-b_O)^T\lambda(t_k)>d_{min}
\end{equation}

\indent Since the robot configuration is always changing, there is a cost function to minimize the bounding box. The constraint can be thought of as mapping the current translation into the obstacle space and seeing how far away it is from the set.\\
\indent In terms of limitations of this paper's approach, the jumping optimization is also performed offline. This is expected, as this optimization is more complicated than the last paper's optimization. Another limitation of this approach is that jumping through obstacles can only be performed statically. In essence, the robot has to be in a certain starting position and a certain landing position, kind of like if a human jumped in place. An advantage of this approach, that the other papers didn't have, is that this was the only paper to go through obstacles instead of just over them. Also, in terms of jumping performance, the quadruped was able to experimentally jump higher than the feasibility approach from the first paper on optimized jumping.\\
\indent In terms of unanswered questions and future works, the authors want to perform more research on optimizing jumps where the landing and takeoff heights are not the same. Another possible unanswered question would be if this approach is valid for a running jump. As in, would the quadruped still be able to avoid the obstacle if the standing jump constraint was removed?

\subsubsection{An Optimal Motion Planning Framework for Quadruped Jumping}
This paper focuses on a different approach to optimization. While the previous papers focused on gradient-based methods, this paper focuses on heuristic-based methods to overcome the highly nonlinear objective function with greater convergence quality. Similarly to the previous papers, this approach also requires an offline optimization but the authors pre-trained different kinds of jumps for the quadruped to use online. This is a step further than the previous papers where their experiments were in highly controlled environments whereas this pre-training methodology allows for a more robust robot.\\
\indent The optimization formulation is poised as:\\
\begin{align}
    \min_x \ \textbf{Fitness}(x)\\
    x(t_{k+1})=x(t_k)+\Delta t \dot{x}(t_k)\\
    \dot{x_{k+1}} = f(u_k, x_k, p_k)\\
    x_k \in \mathbf{X},k=1,2,...,N\\
    u_k\in \mathbf{U}, k = 1,2,...,N-1\\
    x(t_0)=x_0,x(t_{end})=x_{end}
\end{align}
\indent The constraints are kinematic except for $x_k \in \mathbf{X}$ and $u_k\in \mathbf{U}$ where they are kino-dynamic in nature. Essentially, they relate to physical and obstacle constraints where $\mathbf{X}$ and $\mathbf{U}$ are the feasible sets for these constraints. The Fitness objective function is proportional the number of constraints violated. Given some state, x, that violates many constraints then it would have a relatively high Fitness value. Thus, we seek to minimize this value. \\
\indent The optimization variables are reformulated as ground reaction force (GRF) profiles because the GRF at each foot provides all the information needed to generate a trajectory. Similarly to the first jumping paper presented, these force profiles can be modeled as a time-varying polynomial:\\
\begin{equation}
  f_i=\begin{cases}
    \eta_1 \lambda_1[t\quad 1]^T \quad \quad \ \  t\in [0,T_1]\\
    \eta_2 \lambda_2[t^2 \quad t \quad 1]^T \quad t\in [T_1,T_2]\\
    0 \quad \quad \quad \quad \quad \quad \quad \ t\in [T_2,T_3]
  \end{cases}
\end{equation}

\indent Where $f_i$ is the GRF for each foot. $T_i$ represents the unknown time that each phase of jumping and these are concatenated into a vector $T_p$. For example, from time 0 to $T_1$, this is the phase when the quadruped has both feet on the ground. The actual optimization variables are the coefficients $\eta_i\lambda_i$. Since these have no physical meaning, these are again reformulated in terms of the state and state velocities. In the end, the entire decision vector is:
\begin{equation}
    D_{opt}^* = [T_p \quad x(T_1) \dot{x}(T_3)]\in \mathbf{R}^{12}
\end{equation}

\indent The actual optimization scheme is a differential evolution algorithm. Essentially, we initialize many random guesses as a candidate solution to the objective function. Then, we keep perturbing each of these random guesses and we only keep the best solutions (the candidates that have the lowest objective function). We keep perturbing and keeping the best solutions until convergence is achieved. Put simply, this is a survival of the fittest optimization scheme. More formally:\\
\\
\textit{Randomly initialize population vector};\\
\indent \textbf{for} $g\leftarrow 1$ \textbf{to} M \textbf{do}\\
\indent \indent \textbf{for} $i\leftarrow 1$ \textbf{to} N \textbf{do}\\
\indent \indent \textit{Mutation and Crossover};\\
\indent \indent \indent \textbf{for} $j \leftarrow 1$ \textbf{to} W \textbf{do}\\
\indent \indent\indent \indent $v_{i,j} \leftarrow M(x_{i,j}(g))$;\\
\indent \indent\indent \indent $u_{i,j} \leftarrow C(x_{i,j}(g),v_{i,j}(g))$;\\
\indent \indent\indent \textbf{end}\\
\indent \indent\indent \textit{Selection};\\
\indent \indent\indent \textbf{if} $Fitness(Ui(g),k)<Fitness(Xi(g),k)$\\
\indent \indent\indent \indent \textbf{then}\\
\indent \indent\indent \indent \indent $X_i(g) \leftarrow U_i(g)$;\\
\indent \indent\indent \indent \textbf{if}($Fitness(Xi(g),k)<Fitness(Dopt(g),k)$)\\
\indent \indent\indent \indent \indent \textbf{then} $D_{opt}\leftarrow X_i(g)$;\\
\indent \indent\indent \textbf{else}\\
\indent \indent\indent \indent $X_i(g)\leftarrow X_i(g)$;\\
\indent \indent\indent \textbf{end}\\
\indent \indent \textbf{end}\\
\indent $g\leftarrow g+1$;\\
\textbf{end}\\
\\
\indent In terms of limitations of this paper, it can only perform jumps if it is pre-trained on the jump and added onto its internal library of motion planning. Another limitation, similar to previous papers, is that we still assume a planar version of dynamics. In terms of advantages, compared to the other papers,  the convergence quality is much better because of the heuristic-based optimization approach. With a large variance in initial guesses, the solver can pass over local minima which is an advantage over gradient-based methods. Another advantage over the other papers is that it is robust against window-shaped obstacles and does not use prior knowledge like reference trajectories and contact scheduling.\\
\indent In terms of unanswered questions and future works, this approach does not take landing accuracy into account so some future work can be revolved around generalizing the landing approach. One possible approach to solve this problem is using MPC. By utilizing a predictive model of the robot's dynamics, MPC optimizes a sequence of control actions over a finite time horizon. The process involves system modeling, state estimation, trajectory planning, cost function design, optimization, and real-time implementation. The dynamic model captures the robot's physical properties, while state estimation provides accurate information about the robot's state. Trajectory planning incorporates task-specific requirements and constraints, and a carefully designed cost function quantifies the desired performance. Optimization techniques are used to find the sequence of control actions that minimize the cost function while satisfying system dynamics and constraints. The optimized control actions are applied in real-time, with the MPC algorithm continuously re-planning and adapting to changes. This enables legged robots to achieve safe and stable landings by anticipating future behavior, accounting for uncertainties, and actively adjusting control actions.

\subsubsection{Contact-timing and Trajectory Optimization for 3D Jumping on Quadruped Robots}
This paper focuses on programmatically optimizing the contact scheduling of each jump. This is a further improvement on the previous papers as it is common to manually tune a predefined interval for each jumping phase (e.g. 0.5 seconds for the first phase, 0.6 seconds for the second phase).\\
This contribution is relevant because the manual tuning of contact scheduling is often time-consuming and nonoptimal. Not just that, but optimal contact schedules play a role in reducing actuator effort as well.
\indent The contact timing optimization is formulated as:
\begin{align}
    \min_{x,f,e_R} \sum_{k=1}^N \epsilon_\Omega \Omega_k^T\Omega_k+\epsilon_f f_k^Tf_k^T+\epsilon_R e_{R_k}^Te_{R_k}\\
    [R,p,p_f^s](k=1)=[R_0,p_0,p_{f,0}^s]\\
    [\Omega, \dot{\Omega}, \dot{p}](k=1)=0\\
    [R,p,p_f^s](k=N^c)=[R_g,p_g,p_{f,g}^s]\\
    |R(k)[p_f^s(k)-p(k)]-\bar{p_f^s}(k)|\le r\\
    |f_k^{s,x}/f_k^{s,z}|\le \mu , |f_k^{s,y}/f_k^{s,z}|\le \mu \\
    f_{k,min}^s \le f_k^{s,z}\le f_{k,max}^s \\
    p_{k,min}\le p_k \le p_{k,max}\\
    \gamma(x_k, x_{k+1}, f(k), p_f^s(k))=0\\
    R_{k+1}=R_k\exp(\Omega_k T_i/N_i)\\
    \sum_{i=1}^{n_p}T_i \in [T_{min},T_{max}],N^c=\sum_{i=1}^{n_p}N_i\\
    for \quad k=1,2,...,N^c
\end{align}
\indent Position/velocity, rotation/angular velocity, and contact timings are all optimized in this section. In this least squares optimization, the objective function seeks to minimize the angular velocity, ground reaction force for each leg, and the error between the reference rotation trajectory. The optimized ground reaction forces are able to generate a trajectory of the center of mass so the output of this optimization will be a reference trajectory for the jump as well as contact timings for each phase of the jump.\\
\indent These are then fed into the trajectory optimization where the torques are optimized for the least effort. The objective function is:
\begin{multline}
J=\sum_{h=1}^{N-1} \epsilon_q(q_h-q_{ref})^T(q_h-q_{ref})+\\ 
    \epsilon_\tau \tau_h^T\tau_h+\epsilon_N(q_N-q_N^d)^T(q_N-q_N^d)    
\end{multline}

with constraints:
\begin{itemize}
    \item Full body dynamics constraints
    \item Initial configuration: $q(h=0)=q_0,\dot{q}(h=0)=0$
    \item Pre-landing configuration: $q_{j,h}=q_{h,N}^d,\dot{q}_{j,h}=0$
    \item Final Configuration: $q_h(h=N)=q_N^d$
    \item Joint angle constraints: $q_{j,min} \le q_{j,h} \le q_{j,max}$
    \item Joint velocity constraints: $|\dot{q}_{j,h}|\le \dot{q}_{j,max}$
    \item Joint torque constraints : $|\tau_{max}|\le \tau_{max}$
    \item Friction cone limits: $|F_h^x/F^z_h|\le \mu,|F_h^y/F^z_h|\le \mu$
    \item Minimum GRF: $F_h^z\ge F_{min}^z$
\end{itemize}

\indent This formulation is similar to the contact timings optimization, it's a constrained weighted least squares problem to solve for optimized torques. Since this is a torque optimization, the constraints are centered around the hardware limitations like joint angles and initial/final configurations of the quadruped. Once these torques are solved, then these reference profiles are fed into the controller to interface with the hardware.\\
\indent In terms of limitations of this approach, contact timing optimization on the full nonlinear dynamics is a highly complex problem. In the authors' experience, their implementation does not produce a feasible solution for many possible 3D jumps. Because of this, the authors opted to use simplified dynamics to reduce complexity. In terms of advantages, the contact timing is now automated and no longer needs to be manually tuned. This is a huge advantage over previous papers where a key part of the process is finding the correct contact scheduling that works with specific jumps.\\
In terms of unanswered questions and future works, the authors hope to implement a vision to autonomously detect and avoid obstacles in this optimized contact timings framework.

\subsubsection{Future Directions for Gait Planning}
Some future directions, in general, for the subfield of gait planning include implementing a vision for autonomous gait planning and the online viability of these algorithms. Right now, most of these algorithms have to be optimized offline and ported into a controlled testing environment. Only for the most simple objective functions, where it is just a feasibility problem, is where online computation remotely feasible. Even in this edge case, the optimization is fast but the performance is suboptimal. I believe that a problem formulation or hardware innovation that enables high performance while optimizing online is the future of this field as it is a quality that every paper lacks so far.\\
\indent Another path, coupled with online optimization, is implementing vision along with these optimization algorithms. This is probably a hurdle that many papers have not addressed because it is already a large computational cost to implement software that provides an intermediate data representation of video footage that can be passed along to the backend optimization. If these two problems of vision and online optimization can be both solved, then truly robust robots will be created.

\section{Societal Impact}
\indent The potential societal impact of humanoid robots is vast. Though still in the early stages of development, these robots can assist our society with a wide array of tasks. These can be categorized into situations where there is danger to humans, or where the task requires lots of repetitive action. For example, search and rescue during natural disasters or the collapse of mines requires humans to venture into dangerous environments, as does dealing with leaks or catastrophic failures in a nuclear power plant. In these situations, humanoid robots can provide a net positive impact. For cases where repetitive and simple tasks need to be performed, for example in factories or on farms, robots can be used to increase productivity as there will be less need to hire people to do those tasks. However, this leads to potential negative impacts on society such as the replacement of humans and therefore job losses. These factors will have to be carefully weighed by society as the capabilities of humanoid robots expand to a point where they can be reliably deployed in real-world settings.  

\section{Possible solutions to the unanswered questions}
Convex optimization offers potential solutions for the unresolved questions around the stability and control of legged robots. As a mathematical tool, convex optimization is well-suited to handle uncertainties, disturbances, and efficiency in control systems. For instance, we can formulate robust control strategies as convex optimization problems that seek to minimize the impact of uncertainties and disturbances. Furthermore, the computational efficiency of optimization algorithms can enhance real-time control performance. For dealing with uneven terrains, a convex problem can be framed to minimize the deviation of the robot's foot placements from optimal values. Energy efficiency in control systems can also be targeted by formulating a cost function that considers both stability and energy usage. Convex optimization can help integrate learning-based methods into traditional control strategies by defining suitable cost functions and constraints that ensure system stability while allowing for adaptation. Slippage can be handled by formulating constraints in the optimization problem that limit the feasible region of the foot forces. Similarly, switching between different locomotion modes and transitioning between terrains can be managed by defining mode-specific and terrain-specific constraints and cost functions. Force feedback control can be incorporated by adding additional constraints to the optimization problem to ensure balance and stability. Finally, convex optimization allows us to explore various design and control parameters systematically, balancing stability, performance, and efficiency for specific tasks or environments.

\vspace{12pt}

\end{document}